\documentclass[11pt,a4paper]{amsart}
\usepackage{appendix}
\usepackage{amsmath,amsfonts,amsthm,amssymb,amscd}
\usepackage[latin1]{inputenc}
\usepackage{mathrsfs}
\usepackage{latexsym}
\usepackage{graphicx}
\usepackage{enumerate}
\def\Z{\mathbb{Z}}
\def\M{\mathbb{N}_{0}}
\def\N{\mathbb{N}}
\newcommand{\mex}{\operatorname{mex}} 
\begin{document}
\title{2-Pile Nim With a Restricted number of Move-size Imitations}

\author{Urban Larsson}

\email{urban.larsson@chalmers.se, hegarty@chalmers.se} 
\address{Mathematical Sciences,
Chalmers University of Technology and University of 
Gothenburg, G\"oteborg, Sweden}
\date{\today }

\begin{abstract}
We study a variation of the combinatorial game of 2-pile Nim. Move as
in 2-pile Nim but with the following constraint:

Suppose the previous player has just removed say $x>0$ tokens from the 
shorter pile (either pile in case they have the same height). If the next 
player now removes $x$ tokens from the larger pile, then he imitates 
his opponent. For a predetermined natural number $p$, by the rules of the game, 
neither player is allowed to imitate his opponent 
on more than $p-1$ consecutive moves.

We prove that the strategy of this game resembles closely that of a variant of
Wythoff Nim---a variant with a blocking manoeuvre on $p-1$ diagonal positions. 
In fact, we show a slightly more
general result in which we have relaxed the notion of what an imitation 
is.
\end{abstract}
\maketitle
\vspace{-0.5cm}
\begin{center}
\small\textsc{with an appendix by Peter Hegarty}
\end{center}
\vskip 30pt
\section{Introduction} 

A finite impartial game is usually a game where
\begin{itemize}
\item there are 2 players and a starting position,
\item there is a finite set of possible positions of the game,
\item there is no hidden information,
\item there is no chance-device affecting how the players move,
\item the players move alternately and obey the same game rules,
\item there is at least one final position, from which a player
cannot move, which determines the winner of the game and
\item the game ends in a finite number of moves, no matter how it is played.
\end{itemize}
If the winner of the
game is the player who makes the final move, then we play under normal play
rules, otherwise we play a mis\`ere version of the game.

In this paper \emph{a game}, say $G$, is always a
finite impartial game played under normal rules. The player who made
the most recent move will be denoted by \emph{the previous player}. 
A position from which the previous player 
will win, given best play, is called a \emph{$P$-position}, or just $P$. 
A position from which the \emph{next player} will win is called 
an \emph{$N$-position}, or just $N$. The set of all $P$-positions will 
be denoted by $\mathcal{P}=\mathcal{P}_{G}$ and the set of 
all $N$-positions by $\mathcal{N}=\mathcal{N}_{G}$.

Suppose $A$ and $B$ are the two piles of a 2-pile take-away game,
which contain $a\ge 0$ and $b\ge 0$ tokens respectively.
Then the \emph{position} is $(a,b)$ and a \emph{move} 
(or an \emph{option}) is denoted by $(a, b) \rightarrow (c, d)$, where  
$ a - c \ge 0$ and $b - d \ge 0$ but not both $a = c$ and $b = d$.  
 All our games are symmetric in the sense that 
 $(a,b)$ is $P$ if and only if $(b,a)$ is $P$. Hence, 
to simplify notation, 
when we say $(a,b)$ is $P$ ($N$) we also mean $(b,a)$ is $P$ ($N$).
Througout this paper, we 
let $\M$ denote the non-negative integers and $\N$ the 
positive integers. For integers $a<b$ we let $[a,b]$ 
denote the set $\{a,a+1,\ldots ,b\}$.

\subsection{The game of Nim} 
The classical game of Nim is played on a positive number of piles,
each containing a non-negative number of tokens, where 
the players alternately remove 
tokens from precisely one of the non-empty piles---that is, at 
least one token and at most the entire pile---until all piles are gone. 
The winning strategy of Nim is, whenever possible, to move so that 
the ``Nim-sum'' of the pile-heights equals zero, 
see for example \cite{Bou} or \cite{SmSt} (page 3).
When played on one single pile there are only next player winning
positions except when the pile is empty. When played on two
piles, the pile-heights should be equal to ensure victory for the 
previous player.  

\subsection{Adjoin the $P$-positions as moves}
A possible extension of a game is $(\star )$ \emph{to adjoin the
$P$-positions of the original game as moves in the new game}. Clearly
this will alter the $P$-positions of the original game. 

Indeed, if we adjoin 
the $P$-positions of 2-pile Nim as moves, then we get 
another famous game, namely Wythoff Nim (a.k.a Corner the queen), 
see \cite{Wy}. The set 
of moves are: Remove any number of tokens from one of the piles, or remove 
the same number of tokens from both piles. 

The $P$-positions of this game are 
well-known. Let $\phi =\frac{1+\sqrt{5}}{2} $ denote the golden ratio. 
Then $(x,y)$ is a $P$-position if and only if 
$$\left(x,y\right) \in \left\{\, 
\left(\lfloor n\phi \rfloor, \lfloor n\phi ^2\rfloor
\right)\, \mid n \in \M \right\}.$$
We will, in a generalised form, return to the nice 
arithmetic properties of this and other sequences 
in Proposition 1 (see also \cite{HeLa} for further generalisations).

Other examples of $(\star )$ 
are the Wythoff-extensions of $n$-pile Nim for $n\ge 3$ 
discussed in \cite{BlFr, FrKr, Su, SuZe} as well as 
some extensions to the game of 2-pile Wythoff Nim in \cite{FraOz}, 
where the authors adjoin subsets of the Wythoff Nim $P$-positions as moves 
in new games.

\subsection{Remove a game's winning strategy}

There are other ways to construct 
interesting extensions to Nim on just one or two
piles, for example we may introduce a so called \emph{move-size
dynamic} restriction, where the options in some specific way depend on how the
previous player moved (for example how many tokens he removed), or
``pile-size dynamic''\footnote{We understand that pile-size
  dynamic games are not `truly' dynamic since for any given position 
of a game, one may classify each $P$-position without any knowledge 
of how the game has been played up to this point. } restrictions, where 
the options depend on the number of tokens in the respective piles.

 The game of 
``Fibonacci Nim'' in \cite{BeCoGu} (page 483) is a beautiful example of a 
move-size dynamic game on just one pile. This game has been
generalised, for example in 
 \cite{HoReRu}. Treatments of two-pile move-size 
dynamic games can be found in \cite{Co}, extending the (pile-size dynamic) 
``Euclid game'', and in \cite{HoRe}.

The games studied in this paper are move-size dynamic.
In fact, similar to the idea in Section 1.2, 
there is an obvious way to alter the $P$-positions of a game, 
namely $(\star\star)$ \emph{from the original game, 
remove the next-player winning strategy}. For 2-pile Nim this means
that we remove the possibility to \label{imi}
\emph{imitate} the previous player's move, where imitate has the 
following interpretation:\\

\noindent{\bf Definition 0} 
Given two piles, $A$ and $B$, where $\# A\le \# B$---and where the
number of tokens in the respective pile is counted before
the previous player's removal of tokens---then, if the previous
player removed tokens from pile $A$,
the next player \emph{imitates} the previous player's move if he
removes the same number of tokens from pile $B$
as the previous player removed from pile $A$.
\\

This game, we call \emph{Imitation Nim}. 
The intuition is, given the position $(a,b)$, where $a\le b$, Alice
 can prevent Bob from going to $(c,d)$, where $c<a$ and $b-a = d-c$,
 by moving $(a,b)\rightarrow (c,b)$. We illustarate with an example:
\\
\\
{\bf Example 1} 
Suppose the game is Imitation Nim and the
position is $(1,3)$. If this is an a initial position, 
then there is no `dynamic' restriction on the
next move so that the set $\{(1,2),(1,1),(1,0),(0,3)\}$ of Nim 
options is identical to the set of Imitation Nim options. 
But this holds also, if the previous player's move was
\begin{equation}\label{move1}
(1,x)\rightarrow (1,3),
\end{equation}
or
\begin{equation}\label{move2}
(x,3)\rightarrow (1,3),
\end{equation}
where $x\ge 4.$
For these cases, the imitation rule does not apply since the 
previous player removed tokens from the pile with more tokens.

If on the other hand, the previous move was 
\begin{equation}\label{move3}
(x,3)\rightarrow (1,3),
\end{equation} 
 where $x\in\{2,3\}$ then, by the imitation rule, precisely the 
option $(1,3)\rightarrow (1,3-x+1)$ is prohibited. 

Further, $(3,3)\rightarrow (1,3)$ is a 
losing move---since, as we will see in Proposition 0(i), 
$(1,3)\rightarrow (1,2)$ is a winning move. But, by the imitation 
rule, $(2,3)\rightarrow(1,3)$ is a winning 
move---since for this case $(1,3)\rightarrow(1,2)$ is forbidden.\\

This last observation leads us to ask a general question for a move-size
dynamic game, roughly:
\emph{When does the move-size dynamic rule change the outcome of a game?} 
To clarify this question, let us introduce 
some non-standard terminology, valid for any move-size dynamic game.
\\

\noindent {\bf Definition 1}
Let $G$ be a move-size dynamic game. A position $(x,y)\in G$ is 
\begin{enumerate}
\item \emph{dynamic}: if, in the course of 
the game, we cannot tell whether it is $P$ or $N$ without 
knowing the history---at least the most recent move---of the game;
\item \emph{non-dynamic}
\begin{enumerate} 
\item [\emph{$P$}:] if it is  
  $P$ regardless of any previous move(s),
\item [\emph{$N$}:] ditto, but $N$.
\end{enumerate} 
\end{enumerate}
\vspace{2 mm}
 
\noindent{\bf Remark 1}
Henceforth, if not stated otherwise, we will
 think of a (move-size dynamic) 
game as a game where the progress towards the current position is memorized 
in an appropriate manner. A consequence of 
this approach is that each (dynamic) position is $P$ or $N$.
\\

In light of these definitions, we will now characterize the 
winning positions of a game of Imitation Nim (see also Figure 1)---this 
is a special case of our 
main theorem in Section 2, notice for example the absence 
of Wythoff Nim $P$-positions that are dynamic, considered as positions of 
Imitation Nim.
\\

\noindent {\bf Proposition 0} Let $0\le a\le b$ be integers. Suppose 
the game is Imitation Nim. Then $(a, b)$ is
  \begin{enumerate}
  \item[(i)] non-dynamic $P$ if and only if it is a $P$-position 
as of Wythoff Nim;
  \item[(ii)] non-dynamic $N$ if and only if
    \begin{enumerate}
      \item there are integers $0\le c\le d<b$ 
        with $b-a = d-c$ such that $(c,d)$ is a $P$-position 
        of Wythoff Nim, or
      \item there is a $0\le c< a$ such that $(a, c)$ 
        is a $P$-position of Wythoff Nim.
    \end{enumerate}
\end{enumerate}
\vspace{2 mm}

\begin{figure}
\begin{center}
	  \resizebox{10cm}{!}{\includegraphics{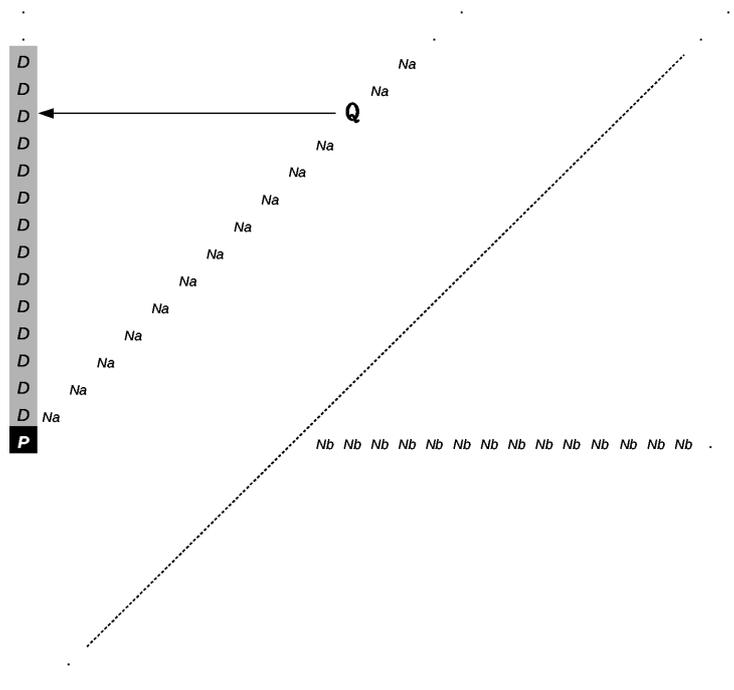}}  
\end{center}\caption{The strategy of Imitation Nim. The P is a (Wythoff Nim) 
$P$-position north of the main diagonal. 
The D:s are dynamic positions. The arrow symbolises a winning move from Q. 
The Na:s are 
the positions of type (iia) in Proposition 0, the Nb:s of type (iib).}
\end{figure}

\noindent {\bf Remark 2} Given the notation in Proposition 0, 
it is wellknown (see also Figure 1) that: 
There is an $x<a$ such that $(x,b)$ 
is a $P$-position of Wythoff Nim implies (iia). One may also note that, 
by symmetry, there is an intersection of type (iia) and (iib) 
positions, namely whenever $a=d$, that is whenever $c<a<b$ is an arithmetic 
progression.\\

By Proposition 0 and Remark 2 $(c,b)$ is a dynamic position of Imitation Nim 
if and only if there is a $P$-position of Wythoff Nim, $(c,d)$, 
with $c\le d<b$. Further, with notation as in (iia), we get 
that $(c,b)$ is dynamic $P$ if and only if 
the previous player moved $(a,b)\rightarrow (c,b)$.
\\

Recall that the first few $P$-positions of Wythoff Nim are $$(0,0),
(1,2), (3,5), \ldots .$$ Hence, in Example 1, a
(non-dynamic) $P$-position of Imitation Nim is $(1,2)$. The position
$(1,3)$ is, by the comment after Remark 2, dynamic. 
The positions $(2,3)$ and $(3,4)$ are non-dynamic $N$-positions,
by Proposition 0(iia). As examples of 
non-dynamic $N$-positions of type (iib), we may take $(2,x)$, $x\ge 3$. 
By the comment after Remark 2 (again), we get:\\

\noindent {\bf Corollary 0} Treated as initial positions, the
$P$-positions of Imitation Nim are identical to the $P$-positions of
Wythoff Nim.
\\

\noindent {\bf Remark 3} For a given position, the rules of Wythoff
Nim allow more options than those of Nim, whereas the rules of Imitation
Nim give fewer. Nevertheless, the $P$-positions are identical 
if one only considers
starting positions. Hence, one might want to view these variants 
of 2-pile Nim as each others ``duals''.

\subsection{Two extensions of Imitation Nim and their ``duals''}
Our first reference for a move-size dynamic game is \cite{BeCoGu}. 
But we have not been able to find any literature on the subject 
of \emph{games with memory}, which is our next topic. 
\subsubsection{A game with memory} 
A natural extension of Imitation Nim is, given
 $p\in \N$, to allow $p - 1$ consecutive imitations---by one and the 
same player---but to prohibit the $p$:th imitation.
We denote this game by \emph{$(p,1)$-Imitation Nim}.
\\

\noindent {\bf Remark 4} This rule removes the winning strategy 
from 2-pile Nim if and only if the number of tokens in each
pile is $\ge p$.
\\

{\noindent \bf Example 2}
Suppose the game is $(2,1)$-Imitation Nim (so that no two
consecutive imitations by one and the same player are 
allowed). Suppose the starting position is $(2,2)$ and that
Alice moves to $(1,2)$. Then, if Bob moves to $(1,1)$,
Alice will move to $(0,1)$, which is $P$ for a game with this 
particular history. This is because the move 
 $(0,1)\rightarrow (0,0)$ would have been 
a second consecutive imitation for Bob and hence is no option.
If Bob chooses instead to move to $(0,2)$
then Alice can win in the next move, since $2>1$ (so the
imitation rule does not apply).

Indeed, Alice's first move is a winning move, 
so $(2,2)$ is $N$ (which is non-dynamic) and $(1,2)$ is $P$. But, 
if $(1,2)$ would have been an initial position then it
would have been $N$, since $(1,2)\rightarrow (1,1)$ would
have been a winning move. So $(1,2)$ is dynamic. 

Clearly $(0,0)$ is non-dynamic $P$. Otherwise the
'least' non-dynamic $P$-position is $(2,3)$, since $(2,2)$ is $N$ and
$(2,1)$ or $(1,3) \rightarrow (1,1)$ would be winning moves, as would 
$(2,0)$ or $(0,3)\rightarrow (0,0)$. 

\subsubsection{The dual of (p,1)-Imitation Nim}
In \cite{HeLa, Lar} we put a \emph{Muller twist} or \emph{blocking manoeuvre} 
on the game of Wythoff Nim. For an introduction to the concept
of a blocking manoeuvre, see 
for example \cite{SmSt}. Variations on Nim with a Muller twist can also be 
found, for example, in \cite{GaSt} (which generalises a result in \cite{SmSt}), 
\cite{HoRe1} and \cite{Zh}.

Fix two positive integers $p$ and $m$. Suppose the
current pile-position is $(a, b)$. The rules are: Before the next player
removes any tokens, the previous player is allowed 
to announce $j \in [1, p-1]$
positions, say $(a_1, b_1), \ldots, (a_j, b_j)$ 
where $b_i - a_i = b - a$, to which the next player may not move. Once the
next player has moved, any blocking manoeuvre is forgotten. 
Otherwise move as in Wythoff Nim.

We will show that as a generalisation of Corollary 0, 
if $X$ is a starting position of $(p,1)$-Imitation Nim then 
it is $P$ if and only if it is a $P$-position of $(p,1)$-Wythoff Nim. 
Further, a generalisation of Proposition 0 holds, but let us now move on 
to our next extension of Imitation Nim.

\subsubsection{A relaxed imitation}
Let $m\in \N$. 
We relax the notion of an imitation to an \emph{$m$-imitation} (or
just imitation) by saying: provided the previous player 
removed $x$ tokens 
from pile $A$, with notation as in Definition 0, then the 
next player $m$-imitates the previous player's move if he removes 
$y\in [x, x + m - 1]$ tokens from 
pile $B$. \\

\noindent {\bf Definition 2}
Let $p\in \N$. We denote by \emph{$(p,m)$-Imitation Nim} 
the game where no $p$ consecutive 
$m$-imitations are allowed by one and the same player.\\

\noindent {\bf Example 3} Suppose that the game is $(1,2)$-Imitation
Nim, so that no 2-imitation is allowed. Then
if the starting position is $(1,2)$ and Alice moves to $(0,2)$, Bob cannot 
move, hence $(1,2)$ is an $N$-position and it must be non-dynamic
since $(1,2)\rightarrow (0,2)$ is always an option regardless of whether
there was a previous move or not.
\subsubsection{The dual of (1,m)-Imitation Nim}
Fix a positive integer $m$. There is a generalisation of Wythoff Nim, 
see \cite{Fra}, here denoted by $(1,m)$-Wythoff Nim, 
which (as we will show in Section 2) has a natural 
$P$-position correspondence with $(1,m)$-Imitation Nim.
 The rules for this game are: remove any number of 
tokens from precisely one of the piles, or remove
tokens from both piles, say $x$ and $y$ tokens respectively, with
the restriction that $|\, x - y\, |  < m$.

 And indeed, to continue Example 3, $(1,2)$ is certainly an
$N$-position of $(1,2)$-Wythoff Nim, since here $(1,2)\rightarrow (0,0)$ is
an option. On the other hand $(1,3)$ is $P$---and 
non-dynamic $P$ of $(1,2)$-Imitation Nim  
 since if Alice moves $(1,3)\rightarrow (0,3)\text{ or }(1,0)$ it does not 
prevent Bob from winning and
$(1,3)\rightarrow (1,2)\text{ or } (1,1)$ are losers, since Bob may take 
advantage of the imitation-rule.

In \cite{Fra}, the author shows that the $P$-positions 
of $(1,m)$-Wythoff Nim are 
so-called ``Beatty pairs'' (view for example the appendix, the
original papers in \cite{Ray, Bea} or in \cite{Fra} (page 355) of
the form $( \lfloor n\alpha\rfloor , \lfloor n\beta \rfloor )$, 
where $\beta = \alpha + m$, $n$ is a non-negative integer and
\begin{equation}
\alpha = \frac{2 - m + \sqrt{m^2 + 4}}{2} . \label{mWN}
\end{equation}
\subsubsection{The $P$-positions of $(p,m)$-Wythoff Nim} 
In the game of $(p, m)$-Wythoff Nim, originally defined 
in \cite{HeLa} (as $p$-blocking $m$-Wythoff Nim), a player may move 
as in $(1,m)$-Wythoff Nim and block positions as in $(p,1)$-Wythoff Nim.
From this point onwards whenever we write Wythoff's game or $W = W_{p,m}$ we 
intend $(p,m)$-Wythoff Nim.

The $P$-positions of this game can easily be calculated by a 
minimal exclusive algorithm (but with exponential complexity 
in succinct input size) as follows: 
Let $X$ be a set of non-negative integers. Define $\mex (X)$ as the least 
non-negative integer not in $X$, formally  
$\mex (X) := \min \{ x\mid x\in\M \setminus X\}$.\\

\noindent{\bf Definition 3} Given positive integers $p$ and $m$, the 
integer sequences $(a_n)$ and $(b_n)$ are: \label{greedy}
\begin{eqnarray*}
a_n & = & \mex \{a_i, b_i\mid 0\le i < n\};\\
b_n & = & a_n + \delta(n),
\end{eqnarray*}
where $\delta (n) = \delta_{p,m} (n) := \left \lfloor \frac
{n}{p}\right \rfloor m$.\\

The next result follows almost immediately from this definition. See 
also \cite{HeLa} (Proposition 3.1 and Remark 1) for further extensions.\\

\noindent {\bf Proposition 1} 
Let $p, m\in \N$.
\begin{itemize}
\item[(a)]  The $P$-positions of $(p,m)$-Wythoff Nim are the pairs 
$( a_i,b_i) $ and $( b_i,a_i) $, $i\in \M$, as in Definition 3;
\item[(b)] The sequences $(a_i)_{i\ge 0}^\infty$ and $(b_i)_{i\ge p}^\infty $ 
partition $\M$ and for $i\in [0,p-1]$, $a_i=b_i=i$;
\item[(c)] Suppose $(a,b)$ and $(c,d)$ are two distinct $P$-positions of 
$(p,m)$-Wythoff Nim with $a \le b$ and $c\le d$. Then $a < c$ implies 
$ b - a \le d - c $ (and $b < d$);
\item[(d)] For each $\delta \in \N $, if $m\! \mid \! \delta $ then 
$\# \{ i\in \M \mid b_i - a_i = \delta \} = p$, 
otherwise $\# \{ i\in \M \mid b_i - a_i = \delta \} = 0$.
\end{itemize}
The $(p,m)$-Wythoff pairs from Proposition 1 may be expressed via Beatty pairs 
if and only if $p\, |\, m$. In that case one can prove via an 
inductive argument that 
the $P$-positions of $(p,m)$-Wythoff Nim are of the form 
$$(pa_n,\, pb_n),\, (pa_n+1,\, pb_n+1),\, \ldots ,\, 
(pa_n+p-1,\, pb_n+p-1),$$ where $(a_n, b_n)$ are the $P$-positions for 
the game $(1,m / p)$-Wythoff Nim (we believe that this fact has not
been recognized elsewhere, at least not in \cite{HeLa} or \cite{Had} 
in its present form).

For any other $p$ and $m$ we did not have a polynomial time algorithm 
for telling whether a given position is $N$ or $P$, until 
recently---while reviewing this article there has 
been progress on this matter,
so there is a polynomial time algorithm, see \cite{Had}. See also
a conjecture in \cite{HeLa}, Section 5, saying in a specific 
sense that the $(p,m)$-Wythoff pairs are ``close to'' the Beatty
pairs $(n\alpha ,n\beta )$ where $\beta=\alpha +\frac{m}{p}$ and
\begin{align*}
\alpha = \frac{2p-m+\sqrt{m^2+4p^2}}{2p},
\end{align*} 
which is settled for the case $m = 1$ in the appendix.
In the general case, as is shown in \cite{Had}, the explicit 
bounds for $a_n$ and $b_n$ 
are $$(n-p+1)\alpha \le a_n \le n\alpha $$ 
and $$(n-p+1)\beta \le b_n \le n\beta .$$ 

A reader who, at this point,  
feels ready to plough into the main idea of our result,  
may move on directly to Section 2---where we state how 
the winning positions of $(p,m)$-Imitation Nim correlate 
to those of $(p,m)$-Wythoff Nim and give a proof for the case $m=1$. 
In Section 3 we finish off with a couple of suggestions for future work.

\subsubsection{Further Examples}
In this section we give two examples of games where $p>1$ and $m>1$ 
(simultaneously), namely in Example 4 $(3,2)$-Imitation Nim 
and in Example 5 $(3,3)$-Imitation Nim. The style is informal.

In Example 4 the winning strategy (via the imitation rule) is in 
a direct analogy to the case $m=1$. In Example 5 we 
indicate how our relaxation of the imitation rule changes how 
a player may take advantage of the imitation rule---in a way 
that is impossible for the case $m=1$. We 
illustrate why this does not affect 
the nice correlation between the winning positions of 
Imitation Nim and Wythoff's game. Hence these 
examples may well be studied in connection with (a second reading of) 
the proof of Theorem 1.\\

\noindent{\bf Example 4} The first 
few $P$-positions of $(3,2)$-Wythoff Nim are
$$(0,0),(1,1),(2,2),(3,5),(4,6),(7,9),(8,12),(10,14),$$
$$(11,15),(13,19),(16,22),(17,23),(18,26),(20,28).$$ For the moment
assume that the first few non-dynamic $P$-positions of
$(3,2)$-Imitation Nim are $(0,0), (3,5), (8,12), (13,19)$ and $(18,26)$.
 
Suppose the position is $(20,27)$. We `suspect' that this is a 
non-dynamic $N$-position
since irrespective of any previous moves, Alice can move 
$(20,27)\rightarrow (17,27)$. This move clearly resets the counter and 
Alice can make sure that Bob will not 
reach the non-dynamic $P$-position $(13,19)$, because then he would
need to imitate Alice's moves 3 times. Is there any other good move
for Bob? Since any other Nim-type move would take him to another
$N$-position (as of Wythoff's game), he must try and rely on the
imitation rule. So he needs to remove tokens from the pile with 17
tokens. But, however he does this, Alice will, by inspection, be able 
to reach a $P$-position (as of Wythoff's game) \emph{without} imitating
Bob. Namely, if Bob moves to $(x,27)$, then Alice next move 
will be $(x,y)$, where $y\le x + 6$ and $10=27-17>6+3=9$, so the move 
is not an imitation. 

By this example we see that the imitation rule is an eminent tool for Alice, 
whereas Bob is the player who 'suffers its consequences'. In the next 
example Bob tries to get around his predicament by hoping that Alice 
would 'rely too strongly' on the imitation rule.\\

\noindent{\bf Example 5} The first 
few $P$-positions of $(3,3)$-Wythoff Nim are
$$(0,0),(1,1),(2,2),(3,6),(4,7),(5,8),(6,12).$$
Suppose, in a game of $(3,3)$-Imitation Nim, the players have moved 
\begin{align*}
\text{Alice}: (6,9)&\rightarrow (5,9)\\
\text{Bob}: (5,9)&\rightarrow (5,6) \text{ an imitation }\\
\text{Alice}: (5,6)&\rightarrow (4,6)\\
\text{Bob}: (4,6)&\rightarrow (3,6) \text{ no imitation}.
\end{align*}

Bob will win, in spite of Alice trying 
to use the imitation rule for her advantage. The mistake is Alice's 
second move, 
where she should change her `original plan' and not continue to try and 
rely on the imitation rule. For the next  
variation Bob tries to 'confuse' Alice's strategy 
by `swapping piles',
\begin{align*}
\text{Alice}: (3,3)&\rightarrow (2,3)\\
\text{Bob}: (2,3)&\rightarrow (2,1).
\end{align*}

Bob has imitated Alice's move once. 
If Alice continues her previous strategy by removing tokens 
from the shorter pile, say by moving $(2,1)\rightarrow (2,0)$, Bob will 
imitate Alice's move a second time and win. Now Alice's correct strategy 
 is rather to remove token(s) from the higher pile,
\begin{align*}
\text{Alice}: (2,1)&\rightarrow (1,1)\\
\text{Bob}: (1,1)&\rightarrow (0,1)\\
\text{Alice}: (0,1)&\rightarrow (0,0).
\end{align*}

Then, Alice has become the player who imitates, but nevertheless wins.

\section{The winning strategy of Imitation Nim}
For the statement of our main theorem we use some more terminology. 
\\

\noindent{\bf Definition 4} Suppose the constants $p$ and $m$ are given as 
in Imitation Nim or in Wythoff's game. Then, if $a,b\in \M $,
\[\xi (a, b)= \xi_{p,m} \bigl((a, b)\bigr) 
:=\#\bigl\{\, (i, j) \in \mathcal{P}_{W_{p,m}} 
\mid \; j - i =b-a,\; i < a\, \bigr\}. \]
\vspace{2 mm}
 
Then according to Proposition 1(d), $$0\le \xi (a,b) \le p, $$  
and indeed, if $(a,b)\in \mathcal{P}_{W_{p,m}} $ then $\xi(a,b)<p$ 
 equals the number of $P$-positions the previous 
player has to block off (given that we are playing Wythoff's game) in 
order to win.\\

\noindent{\bf Definition 5}
Let $(a,b)$ be a position of a game of 
$(p,m)$-Imitation Nim. Put $$L(a,b)=L_{p,m}((a,b)):=p-1$$ if
\begin{itemize}
\item[(A)] $(a,b)$ is the starting position, or 
\item[(B)] $(c,d)\rightarrow (a,b)$ was the most recent move and $(c,d)$ was 
the starting position, or
\item[(C)] The previous move was 
$(e,f)\rightarrow (c,d)$ but the move (or option) 
$(c,d)\rightarrow(a,b)$ is not an $m$-imitation.
\end{itemize}
Otherwise, with notation as in (C), put $$L(a,b)=L(e,f)-1.$$
\vspace{2 mm}
Notice that by the definition of $(p,m)$-Imitation Nim, $$-1\le L(a,b)< p,$$ 
namely it will be convenient to 
allow $L(a,b)=-1$, although a player cannot move 
$(c,d)\rightarrow (a,b)$ if it is an imitation and $L(e,f)=0$.

Indeed $L(e,f)$ represents the number 
of imitations the player moving from $(c,d)$ still has 'in credit'.\\

\noindent{\bf Theorem 1}
 Let $0\le a\le b$ be integers and suppose the game is  $(p,m)$-Imitation Nim. 
Then $(a, b)$ is $P$ if and only if 
 \begin{enumerate} 
\item [(I)] $(a,b)\in \mathcal{P}_{W_{p,m}}$ and $0\le \xi(a,b)\le L(a,b) $, or
\item [(II)] there is a  $a\le c < b$ such that $(a,c)\in \mathcal{P}_{W_{p,m}}$ 
but $-1\le L(a,c)< \xi(a,c)\le p-1$.
\end{enumerate}  

\vspace{2 mm}
\noindent{\bf Corollary 1}
If $(a,b)$ is a starting position of $(p,m)$-Imitation Nim it is $P$ 
if and only if it is a $P$-position of $(p,m)$-Wythoff Nim.\\

\noindent{\bf Proof of Corollary 1} 
Put $L(\cdot)=p-1$ in Theorem 1. $\hfill \Box$ \\

\begin{figure}
\begin{center}
	  \resizebox{10cm}{!}{\includegraphics{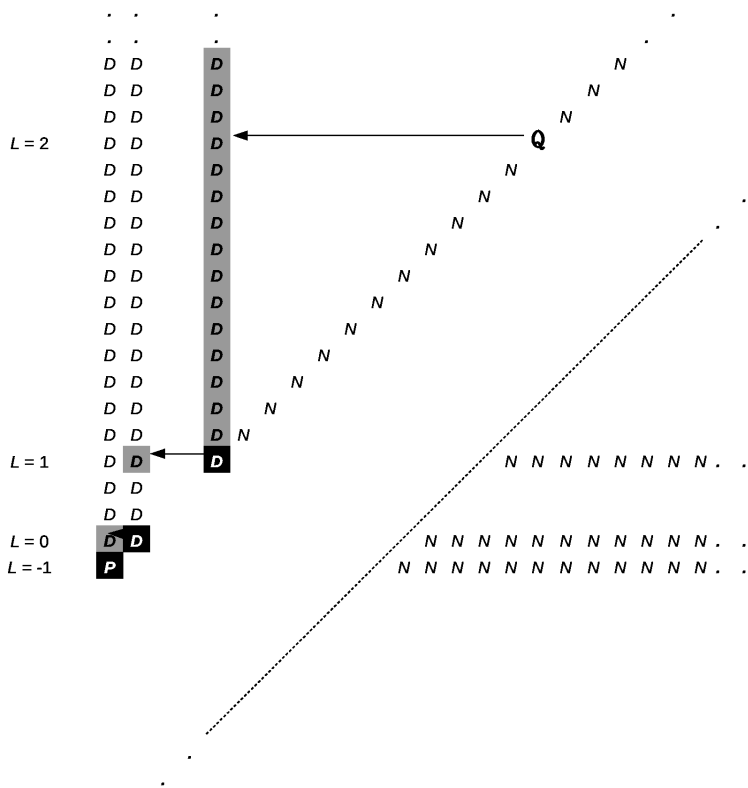}}  
\end{center}\caption{A strategy of a game of $(3,1)$-Imitation Nim. 
The P is a non-dynamic $P$-position north of the main diagonal. The 
black positions are all $P$-positions of $(3,1)$-Wythoff Nim on 
one and the same SW-NE diagonal.
The D:s are dynamic positions. The arrows symbolise 3 consecutive 
winning moves from a position {\bf Q}. 
A position is grey or black if and only if it is $P$ 
in some winning strategy (see also footnote 2). 
The N:s are non-dynamic $N$-positions.}
\end{figure}

By Theorem 1(I) and the remark after Definition 5 we get 
that $(a, b)$ is non-dynamic $P$ if and only if $(a,b)\in \mathcal{P}_{W}$ and 
$\xi(a,b)=0$. On the other hand, if $(a,b)\in \mathcal{N}_W$ it is dynamic 
if and 
only if there is a $a\le c < b$ such that $(a,c)\in \mathcal{P}_{W}$ 
\footnote{While reviewing this article we have found 
out that, under the assumption that through the course of the game at 
least one player has always used a perfect strategy, several dynamic 
positions (only $N$-positions of Wythoff's game though) 'are non-dynamic $N$'. 
 In this sense
 one might want to define the set of \emph{perfect dynamic} positions, 
the subset of dynamic positions that are $P$ in 
some (perfect) strategy, see also Section 3 and Figure 2.}.\\ 

\noindent{\bf Proof of Theorem 1}
We only give the proof for the case $m=1$. In this way 
we may put a stronger emphasis on the idea of the game, at the expense 
of technical details. Whenever we refer to Proposition 1(b, c or d) 
we also intend Proposition 1(a).

Suppose $(a,b)$ is as in (I). Then we need to show that, if $(x,y)$ 
is an option of $(a,b)$ then $(x,y)$ is neither of form (I) nor (II).
 
But Proposition 1(b) gives immediately that $(x,y)\in \mathcal{N}_W$ so 
suppose $(x,y)$ is of form (II). Then there is a $x\le c < y$ 
such that $(x,c)\in\mathcal{P}_W$ and $L(x,c)<\xi(x,c)$. Since, by (I),  
$\xi(a,b)\le L(a,b) $ and $L(a,b)-1\le L(x,c) (\le L(a,b))$ we 
get that $$\xi(a,b) \le L(a,b)\le L(x,c)+1\le\xi(x,c),$$ which, 
in case $c-x=b-a$, is possible 
if and only if $\xi(a,b)=\xi(x,c)$. But then, since, by our assumptions,
$(x,c)\in\mathcal{P}_W$ and $(a,b)\in\mathcal{P}_W$,
we get $(a,b)=(x,c)$, which is impossible. 

So suppose that $c-x\ne b-a$. Then, by Proposition 1(c), $c-x < b-a$. We 
have 2 possibilities:
\begin{itemize}
\item [$y=b$:] Then if $(x,b)\rightarrow (x,c)$ is an imitation 
of $(a,b)\rightarrow (x,b)$ we get $b-c >a-x=b-c$, a contradiction.
\item [$x=a$:] For this case the move $(a,y)\rightarrow (a,c)$ 
cannot be an imitation of $(a,b)\rightarrow (a,y)$ since the 
previous player removed tokens from the larger pile. 
Then $L(a,c)=p-1\ge \xi(a,c)$ since, by (II), $(a,c)\in \mathcal{P}_W$.
\end{itemize}
Hence we may conclude that if $(a,b)$ is of form (I) then an option 
of $(a,b)$ is neither of form (I) nor (II).

Suppose now that $(a,b)$ is of form (II). Then $(a,c)\in \mathcal{P}_W$ 
is an option of $(a,b)$ but we have $L(a,c)<\xi(a,c)$ so $(a,c)$ is not 
of form (I). Since $(a,c)\in \mathcal{P}_W$, by Proposition 1(b), 
it cannot be of form (II).
But then, since $b>c$, by Proposition 1(b) and (c), 
any other option of $(a,b)$, 
say $(x,y)$ must be an $N$-position of Wythoff's game 
so suppose $(x,y)$ is of form (II). We get two cases:
\begin{itemize}
\item [$y=b$:] Then $0\le x<a$ and there is 
an option $(x,d)\in \mathcal{P}_W$ of $(x,b)$ with $x\le d < b$, but 
by Proposition 1(b) and (c) $d-x\le c-a<b-a$ so that 
$(x,b)\rightarrow (x,d)$ does not imitate $(a,b)\rightarrow (x,b)$. 
Hence $L(x,d)=p-1\ge \xi(x,d)$, which contradicts the assumptions in (II).
\item [$x=a$:] Then $0\le y<b$. But then, if $y>c$, 
$(a,c)\in \mathcal{P}_W$ is an option of $(a,y)$ and two consecutive 
moves from the larger pile would give
 $L(a,c)=p-1\ge \xi(a,c)$. Otherwise, by Proposition 1(b), 
there is no option of $(a,y)$ in $\mathcal{P}_W$. In either case a 
contradiction to the assumptions in (II).
\end{itemize}
We are done with the first part of the proof.

Therefore, for the remainder of the proof, assume that $(\alpha,\beta)$, 
$0\le \alpha\le \beta$, is neither of form (I) nor (II). Then
\begin{itemize}
\item[(i)] $(\alpha,\beta)\in \mathcal{P}_{W}$ implies 
$0\le L(\alpha,\beta) < \xi(\alpha,\beta)\le p-1$, and
\item[(ii)] there is a $\alpha\le c < \beta$ such 
that $(\alpha,c)\in \mathcal{P}_{W}$ implies 
$0\le \xi(\alpha,c)\le L(\alpha,c)\le p-1$.
\end{itemize}
We need to find an option of $(\alpha,\beta)$, say $(x,y)$,  
of form (I) or (II). 

If $(\alpha,\beta)\in \mathcal{P}_W$, by Proposition 1(b), 
(ii) is trivially satisfied, and by
(i) $\xi (\alpha ,\beta)>0 $, so there is 
a position $(x,z)\in \mathcal{P}_W$ such that $z-x=\beta-\alpha$ 
with $x\le z< \beta(=y)$. Then, 
since $L(\alpha ,\beta)<\xi(\alpha,\beta)$, 
the option $(x,\beta)$ satisfies (II) 
(and hence, by the imitation rule, $(\alpha,\beta)\rightarrow (x,\beta)$ 
is the desired winning move).

For the case $(\alpha,\beta)\in \mathcal{N}_W$ (here (i) is trivially
true), suppose 
$(\alpha,c)\in \mathcal{P}_W$ with $\alpha\le c < \beta$. 
Then (ii) gives $L(\alpha,c)\ge \xi(\alpha,c)$, which clearly holds 
for example
if the most recent move was no imitation. In any case it imediately
implies (I). 

If $c<\alpha$, with $(\alpha,c)\in \mathcal{P}_W$, then (ii) 
holds trivially by Proposition 1(b) and so (I) holds because
$(\alpha,\beta)\rightarrow (\alpha,c)$ is no imitation (since 
if it was, the previous move must have been from the larger pile).

If $c<\alpha $ with $(c,\beta)\in \mathcal{P}_W$ the move
$(\alpha,\beta)\rightarrow (c,\beta)$ is no imitation
since tokens have been removed from the smaller pile. Hence 
$p-1 = L(c,\beta)\ge \xi(c,\beta)$.

The only remaining case for $(\alpha,\beta )$ an $N$-position of 
Wythoff's game is whenever there is a 
position $(x,z)\in \mathcal{P}_W$ such that 
$x<\alpha$ and 
\begin{align}\label{sd}
\beta-\alpha = z-x.
\end{align}
We may assume there is no $c< \beta$ such 
that $(\alpha,c)\in \mathcal{P}_W$ 
(since we are done with this case). Then (ii) holds trivially and  
by Proposition 1(b) there must be a $c>\beta$ such 
that $(\alpha,c)\in \mathcal{P}_W$. But then,
by Proposition 1(c) and (d), we get $\xi(\alpha,\beta ) = p > 0$ and 
so, since we for this case may take $(x,z)$ such that $p-1 =\xi(x,z)$, we get 
$L(x,z)\le p-2 < \xi(x,z)$. Then, by $(\ref{sd})$, clearly $(x,\beta)=(x,y)$ is 
the desired position of form (II). $\hfill \Box$

\section{Final questions}
Let us finish off with some questions.

\begin{itemize}
\item Consider a slightly different setting of an impartial game, namely where 
the second player does not have perfect information, but the first player 
(who has) is not aware of this fact---similar settings have been discussed 
in for example \cite{BeCoGu, Ow}. We may ask, 
for which games (start with the games we have discussed) is there a 
simple \emph{second player's strategy} which lets him \emph{learn} 
the winning strategy of the game \emph{while playing}---in the sense 
that if he starts a new 'partie' of the same game 
at least 'one move after' the first one, he will win? 
\item Is there a generalisation of Wythoff Nim 
to $n> 2$ piles of tokens (see for example \cite{BlFr, FrKr, Su, SuZe}), 
together with 
a generalisation of 2-pile Imitation Nim, such that the $P$-positions 
correlate (at least as starting positions)?
\item Are there other impartial (or partizan) 
games where an imitation rule corresponds in a natural way 
to a blocking manoeuvre? 
\item Can one formulate a general rule 
as to when such correspondences can be found and when not?
\end{itemize}

\noindent{\bf Acknowledgements} At the Integers 
2007 Conference when I introduced Imitation Nim to Aviezri Fraenkel 
I believe he quickly responded ``\emph{Limitation Nim}''.  
As this name emphasises another important aspect of the game, 
I would like to propose it for the ``dual'' of Wythoff's 
original game, that is whenever one wants to emphasise that 
no imitation is allowed. I would also like to thank A. Fraenkel 
for contributing with some references. 

I would like to thank Peter Hegarty for the comments 
and suggestions he has given during the writing of this paper.

\newpage
\section*{Appendix}
\begin{center}
\small\textsc{Peter Hegarty}
\end{center}

$\;$ \\
The purpose of this appendix is to provide a proof of 
Conjecture 5.1 of \cite{HeLa} in the
case $m = 1$, which is the most natural case to consider. Notation
concerning $\lq$multisets' and $\lq$greedy permutations' is consistent 
with Section 2 of \cite{HeLa}. We begin by 
recalling 
\\
\\
{\sc Definition} : Let $r,s$ be positive irrational numbers with $r < s$. 
Then $(r,s)$ is
said to 
be a {\em Beatty pair} if 
\begin{eqnarray}
{1 \over r} + {1 \over s} = 1.
\end{eqnarray}
{\bf Theorem} {\em Let $(r,s)$ be a Beatty pair. Then the map $\tau :
{\N} \rightarrow {\N}$ given by 
\begin{eqnarray*}
\tau( [nr] ) = [ns], \;\; \forall \; n \in {\N}, \;\;\;\;\;\; 
\tau = \tau^{-1},
\end{eqnarray*}
is a well-defined involution of $\N $. If $M$ is the multiset of 
differences $\pm \{ [ns] - [nr] : n \in {\N} \}$, then 
$\tau = \pi_{g}^{M}$. $M$ has asymptotic density equal to $(s-r)^{-1}$.}
\\
\\
{\sc Proof} : That $\tau$ is a well-defined permutation of $\N $ 
is Beatty's theorem. The second and third assertions are then obvious. 
\\
\\
{\bf Proposition} {\em Let $r < s$ be positive real numbers satisfying (7),
and let $d := (s-r)^{-1}$. Then the following are equivalent
\par (i) $r$ is rational
\par (ii) $s$ is rational
\par (iii) $d$ is rational of the form $\frac{mn}{m^{2}-n^{2}}$ for some 
positive rational $m,n$ with $m > n$.}
\\
\\
{\sc Proof} : Straightforward algebra exercise. 
\\
\\
{\sc Notation} : Let $(r,s)$ be a Beatty pair, $d := (s-r)^{-1}$. We denote 
by $M_{d}$ the multisubset of $\N$ consisting of all differences 
$[ns] - [nr]$, for $n \in {\N}$. We denote $\tau_{d} := 
\pi_{g}^{\pm M_{d}}$.
\par As usual, for any positive integers $m$ and $p$, we denote by 
$\mathcal{M}_{m,p}$ the multisubset of $\Z$ consisting of $p$ copies 
of each 
multiple of $m$ and $\pi_{m,p} := \pi_{g}^{\mathcal{M}_{m,p}}$.  We now 
denote by $M_{m,p}$ the submultiset consisting
of all the positive integers in $\mathcal{M}_{m,p}$ and 
$\overline{\pi}_{m,p} :=\pi_{g}^{\pm M_{m,p}}$. Thus 
\begin{eqnarray}
\overline{\pi}_{m,p} (n) + p = \pi_{m,p}(n+p) \;\;\;\; {\hbox{for all $n \in 
{\N}$}}.
\end{eqnarray}   
Since $\mathcal{M}_{m,p}$ 
has density $p/m$, there is obviously a close relation 
between $M_{m,p}$ and $M_{p/m}$, and thus between the permutations 
$\pi_{m,p}$ and $\tau_{p/m}$. The precise nature of this relationship is, 
however, a lot less obvious on the level of permutations. It is the 
purpose of the present note to explore this matter. 
\\
\\
We henceforth assume that $m=1$. 
\\
\\
To simplify notation we fix a value of $p$. We 
set $\pi := \overline{\pi}_{1,p}$, Note that  
\begin{eqnarray*}
r = r_{p} = {(2p-1) + \sqrt{4p^{2}+1} \over 2p}, \;\;\;\; s = s_{p} = r_{p} + 
\frac{1}{p} = {(2p+1) + \sqrt{4p^{2}+1} \over 2p}.
\end{eqnarray*}
{\sc Further notation} : If $X$ is an infinite multisubset of $\N$ we 
write $X = (x_{k})$ to denote the elements of $X$ listed in 
increasing order, thus strictly increasing order when $X$ is an ordinary 
subset of $\N$. The following four subsets of $\N$ will be of 
special interest :
\begin{eqnarray*}
A_{\pi} := \{ n : \pi(n) > n \} := (a_{k}), \\
B_{\pi} := {\N} \backslash A_{\pi} := (b_{k}), \\
A_{\tau} := \{ n : \tau(n) > n \} := (a^{*}_{k}), \\
B_{\tau} := {\N} \backslash A_{\tau} := (b^{*}_{k}).
\end{eqnarray*}
Note that $b_{k} = \pi(a_{k})$, $b^{*}_{k} = \tau(a^{*}_{k})$ for all $k$. We
set 
\begin{eqnarray*}
\epsilon_{k} := (b_{k}-a_{k}) - (b^{*}_{k} - a^{*}_{k}) = (b_{k}-b^{*}_{k})
- (a_{k}-a^{*}_{k}).
\end{eqnarray*}
{\bf Lemma 1} {\em (i) For every $n > 0$, 
\begin{eqnarray*}
| M_{p} \cap [1,n] | = | M_{1,p} \cap [1,n]| + \epsilon,
\end{eqnarray*}
where $\epsilon \in \{0,1,...,p-1\}$. 
\\ (ii) $\epsilon_{k} \in \{0,1\}$ for all $k$ and if $\epsilon_{k} = 1$ then
$k \not\equiv 0 \; ({\hbox{mod $p$}})$.
\\
(iii) $a^{*}_{k+1} - a^{*}_{k} \in \{1,2\}$ for all $k > 0$ and cannot equal 
one for any two consecutive values of $k$.
\\
(iv) $b^{*}_{k+1} - b^{*}_{k} \in \{2,3\}$ for all $k > 0$.}
\\
\\
{\sc Proof} : (i) and (ii) are easy consequences of the various definitions. 
(iii) follows from the fact that $r_{p} \in (3/2,2)$
and (iv) from the fact that $s_{p} \in (2,3)$.   
\\
\\
{\bf Main Theorem} {\em For all $k > 0$, $|a_{k} - a^{*}_{k}| \leq p-1$.}
\\
\\
{\sc Remark} : We suspect, but have not yet been able to prove, that $p-1$ is 
best-possible in this theorem.
\\
\\
{\sc Proof of Theorem} : The proof is an induction on $k$, which is most easily
phrased as an argument by contradiction. Note that $a_{1} = a^{*}_{1} = 1$. 
Suppose the theorem is false and consider the smallest $k$ for which 
$|a^{*}_{k} - a_{k}| \geq p$. Thus $k > 1$. 
\\
\\
{\em Case I} : $a_{k} - a^{*}_{k} \geq p$. 
\\
\\
Let $a_{k}-a^{*}_{k} := p^{\prime} \geq p$. 
Let $b_{l}$ be the largest element of $B_{\pi}$ in $[1,a_{k})$. Then 
$b^{*}_{l-p^{\prime}+1} > a^{*}_{k}$ and Lemma 1(iv) implies that 
$b^{*}_{l} - b_{l} \geq p^{\prime}$. But Lemma 1(ii) then implies that also
$a^{*}_{l} - a_{l} \geq p^{\prime} \geq p$. Since obviously $l < k$, this
contradicts the minimality of $k$.
\\
\\
{\em Case II} : $a^{*}_{k} - a_{k} \geq p$. 
\\
\\
Let $a^{*}_{k} - a_{k} := p^{\prime} \geq p$. Let $b^{*}_{l}$ be the largest 
element of $B_{\tau}$ in $[1,a^{*}_{k})$. Then $b_{l-p^{\prime}+1} > a_{k}$.
Lemma 1(iv) implies that $b_{l-p^{\prime}+1} - b^{*}_{l-p^{\prime}+1} \geq 
p^{\prime}$ and then Lemma 1(ii) implies that $a_{l-p^{\prime}+1}
- a^{*}_{l-p^{\prime}+1} \geq p^{\prime} - 1$. The only way we can avoid a 
contradiction already to the minimality of $k$ is if all of the following 
hold :
\\
\\
(a) $p^{\prime}=p$.
\\
(b) $b^{*}_{i} - b^{*}_{i-1} = 2$ for $i = l,l-1,...,l-p+2$.
\\
(c) $l \not\equiv -1 \; ({\hbox{mod $p$}})$ and $\epsilon_{l-p+1} = 1$. 
\\
\\
To simplify notation a little, set $j := l-p+1$. Now $\epsilon_{j} = 1$ but
parts (i) and (ii) of Lemma 1 imply that we must have $\epsilon_{j+t} = 0$
for some $t \in \{1,...,p-1\}$. Choose the smallest $t$ for which
$\epsilon_{j+t} = 0$. Thus 
\begin{eqnarray*}
b^{*}_{j} - a^{*}_{j} = b^{*}_{j+1} - a^{*}_{j+1} = \cdots = 
b^{*}_{j+t-1} - a^{*}_{j+t-1} = (b^{*}_{j+t} - a^{*}_{j+t}) - 1.
\end{eqnarray*}
From (b) it follows that 
\begin{equation}
a^{*}_{j+t} - a^{*}_{j+t-1} = 1, \;\; a^{*}_{j+\xi} - a^{*}_{j+\xi-1} = 2, \;
\xi = 1,...,t-1.
\end{equation}
Let $b^{*}_{r}$ be the largest element of $B_{\tau}$ in $[1,a^{*}_{j})$. Then
from (9) it follows that 
\begin{equation}
b^{*}_{r+t} - b^{*}_{r+t-1} = 3, \;\; 
b^{*}_{r+\xi} - b^{*}_{r+\xi-1} = 2, \; \xi = 2,...,t-1.
\end{equation}
Together with Lemma 1(iv) this implies that 
\begin{equation}
b^{*}_{r+p-1} - b^{*}_{r+1} \geq 2p-3.
\end{equation} 
But since $a^{*}_{j} = a_{j} - (p-1)$ we have that $b_{r+p-1} <
a_{j}$. 
Together with (11) this enforces 
$b^{*}_{r+p-1} - b_{r+p-1} \geq p$, and then by Lemma
1(ii) we also have $a^{*}_{r+p-1} - a_{r+p-1} \geq p$. Since it is easily
checked that $r+p-1 < k$, we again have a contradiction to the minimality of
$k$, and the proof of the theorem is complete.    
\\
\\
This theorem implies Conjecture 5.1 of \cite{HeLa}. 
Recall that the $P$-positions
of $(p,1)$-Wythoff Nim are the pairs $(n-1,\pi_{1,p}(n)-1)$ for $n \geq 1$.
\\
\\
{\bf Corollary} {\em With
\begin{eqnarray*}
L = L_p = {s_p \over r_p} = {1 + \sqrt{4p^2 + 1} \over 2p}, \;\;\;
l = l_p = {1 \over L_p} ,
\end{eqnarray*}
we have that, for every $n \geq 1$, 
\begin{eqnarray}
\pi_{1,p}(n) \in \left\{ \lfloor nL \rfloor + \epsilon, 
\lfloor nl \rfloor + \epsilon : \epsilon \in \{-1,0,1,2\} \right\}.
\end{eqnarray}}
{\sc Proof} : We have $\pi_{1,p}(n) = n$ for $n = 1,...,p$, and one checks that
(12) thus holds for these $n$. For $n > p$ we have by (8) that 
\begin{eqnarray}
\pi_{1,p}(n) = \pi(n-p) + p, 
\end{eqnarray}
where $\pi = \overline{\pi}_{1,p}$. 
There are two cases to consider, according as
to whether $n-p \in A_{\pi}$ or $B_{\pi}$. We will show in the former case
that $\pi_{1,p}(n) = \lfloor nL \rfloor + \epsilon$ for some $\epsilon \in
\{-1,0,1,2\}$. The proof in the latter case is similar and will be omitted.
\\
\\
So suppose $n-p \in A_{\pi}$, say $n-p = a_k$. Then 
\begin{eqnarray}
\pi(a_k) = b_k = a_k + (b_k^{*} - a_k^{*}) + \epsilon_k.
\end{eqnarray}
Moreover $a_k^{*} = \lfloor kr_p \rfloor$ and $b_k^{*} = \lfloor ks_p \rfloor$, 
from which it is easy to check that 
\begin{eqnarray*}
b_k^{*} = a_k^{*} L + \delta, \;\;\; {\hbox{where $\delta \in (-1,1)$}}.
\end{eqnarray*}
Substituting into (14) and rewriting slightly, we find that 
\begin{eqnarray*}
\pi(a_k) = a_k L + (a_k^{*}-a_k)(L-1) + \delta + \epsilon_k,
\end{eqnarray*}
and hence by (13) that $\pi_{1,p}(n) = nL + \gamma$ where 
\begin{eqnarray*}
\gamma = (a_k^{*}-a_k - p)(L-1) + \delta + \epsilon_k.
\end{eqnarray*}
By Lemma 1, $\epsilon_k \in \{0,1\}$. By the Main Theorem, 
$|a_k^{*} - a_k| \leq p-1$. It is easy to check that $(2p-1)(L-1) < 1$. Hence
$\gamma \in (-2,2)$, from which it follows immediately that
$\pi_{1,p}(n) - \lfloor nL \rfloor \in \{-1,0,1,2\}$. This completes the proof.
\\
\\
{\sc Remark} : As stated in Section 5 of \cite{HeLa}, computer calculations seem 
to suggest that, in fact, (12) holds with just $\epsilon \in \{0,1\}$. So once
again, the results presented here may be possible to improve upon.\\

\end{document}